\documentclass[12pt]{amsart}
\usepackage{amssymb,amscd}
\usepackage[russian]{babel}
\usepackage{graphicx}
\usepackage{pdfpages}
\usepackage{placeins}
\newtheorem{theorem}{Теорема}[section]

\newtheorem{lemma}{Лемма}[section]

\theoremstyle{definition}
\newtheorem{definition}{Определение}[section]

\theoremstyle{remark}

\ifx\pdfpagewidth\undefined \else
\fi
 \evensidemargin \oddsidemargin

\advance\hoffset by -1in
\advance\voffset by -1in

\usepackage[cp1251]{inputenc}

\author{И.~А.~Иванов-Погодаев, А.~Я.~Канель-Белов}

\thanks{Moscow Institute of Physics and Technology, Bar-Ilan University, Мехмат МГУ}

\address{Московский Физико-Технический Институт}
\email{ivanov.pogodaev@mail.ru}
\address{Bar-Ilan University, Israel, Мехмат МГУ}
\email{kanel@mccme.ru}

\title{Полугруппа путей на семействе равномерно эллиптических комплексов}

\begin{document}

\let \mathbf=\texttt
\tabcolsep 2pt

\begin{abstract}

Это третья часть цикла работ, посвященного конструкции конечно определенной бесконечной нильполугруппы, удовлетворяющей тождеству $x^9=0$. Эта конструкция отвечает на проблему Л.~Н.~Шеврина и М.~В.~Сапира, поставленной, например, в Свердловской тетради \cite{Sverdlovsk}.

Полугруппа реализуется как множество кодировок путей на семействе специальных равномерно-эллиптических комплексов.  


В первой работе цикла <<Конечно определенная ниль полугруппа: комплексы с равномерной эллиптичностью>> (\cite{complex}) была построена последовательность комплексов, обладающих набором геометрических свойств. 

Во второй работе цикла <<Детерминированная раскраска семейства комплексов>> (\cite{coloring}) на вершинах и ребрах построенных комплексов была введена конечная кодировка буквами. Было доказано свойство детерминированности такой раскраски, позволяющего ввести конечное множество определяющих соотношений на множестве слов-кодировок путей на комплексах.



В настоящей работе описывается алгоритм приведения произвольного полугруппового слова к каноническому виду. Также доказывается, что слово, содержащее подслово с периодом $9$ может быть приведено к нулю с помощью определяющих соотношений. Кодировки слов, соответствующие достаточно длинным путям не приводятся к нулю не меняют длины, то есть введенная полугруппа является бесконечной.


Работа была проведена с помощью Российского Научного Фонда, Грант 22-11-00177. Первый автор является победителем конкурса ``Молодая математика России''.








\end{abstract}

\maketitle

Ключевые слова: конечно определенные полугруппы, нильполугруппы, конечно определенные кольца, конечно определенные группы. 

УДК: 512.53 MSC: 20M05

\medskip
\tableofcontents

\medskip

\section{Введение} \label{nachalo}

Это завершающая часть цикла работ посвященного построению конечно определенных нильполугрупп. Доказана следующая

\medskip
{\bf Теорема.} {\it Существует конечно определенная бесконечная нильполугруппа, удовлетворяющая тождеству $x^9=0$.}
\medskip

Ниже приводится краткое содержание первой и второй работы цикла. Введение ко всему циклу дано в первой работе \cite{complex}. Там же приведены ссылки на публикации о родственных задачах и подходах.


\section{Содержание первой работы}

В первой части цикла работ <<Конечно определенная ниль полугруппа: комплексы с равномерной эллиптичностью>> (\cite{complex}) была изложена геометрическая структура семейства комплексов, служащих базой для введения определяющих соотношений. 

Ниже кратко излагаются основные положения первой части работы. Точные определения и доказательства можно найти в \cite{complex}.

\begin{figure}[hbtp]
\centering
\includegraphics[width=0.5\textwidth]{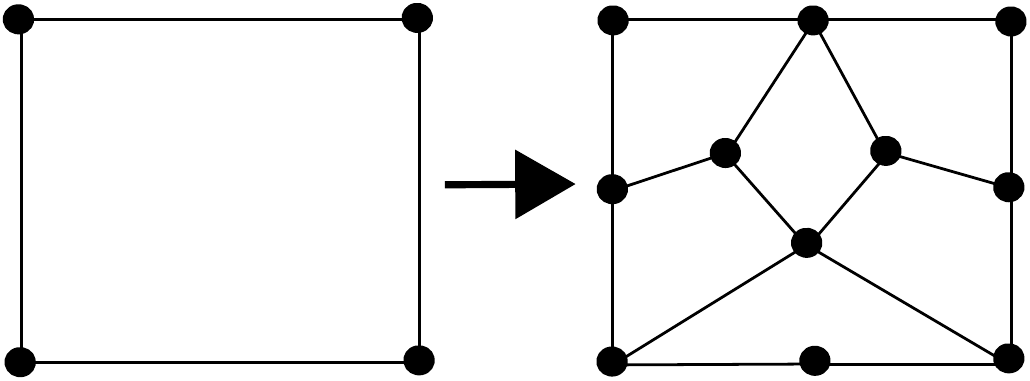}
\caption{Разбиение.}
\label{subst}
\end{figure}

Основной конструкцией работы является семейство геометрических комплексов, каждый из которых состоит из конечного числа квадратов, приклееных друг к другу сторона к стороне специальным образом. Нас будут интересовать свойства комплексов как графов. Комплекс $1$ уровня это просто квадрат из $4$ вершин (плитка). Комплекс $2$ и $3$ уровня изображены на рисунке. Комплекс $3$ уровня получается разбиением каждого квадрата на $6$ квадратов, при этом в середине каждого ребра появляется дополнительная вершина (рисунок~\ref{subst}. 
{\it Макроплиткой $n$ уровня} называется плоская часть комплекса -- результат применения 
$n-1$ уровней разбиения к квадрату из $4$ вершин. 

 Можно заметить, что при этом образуются различные виды вершин --  угловые, краевые, боковые, и внутренние. Угловые и краевые лежат на краю комплекса (в углах, или на сторонах), боковые лежат на ребре, проведенном при некотором разбиении, внутренние -- образуются при некотором разбиении.

\begin{figure}[hbtp]
\centering
\includegraphics[width=0.5\textwidth]{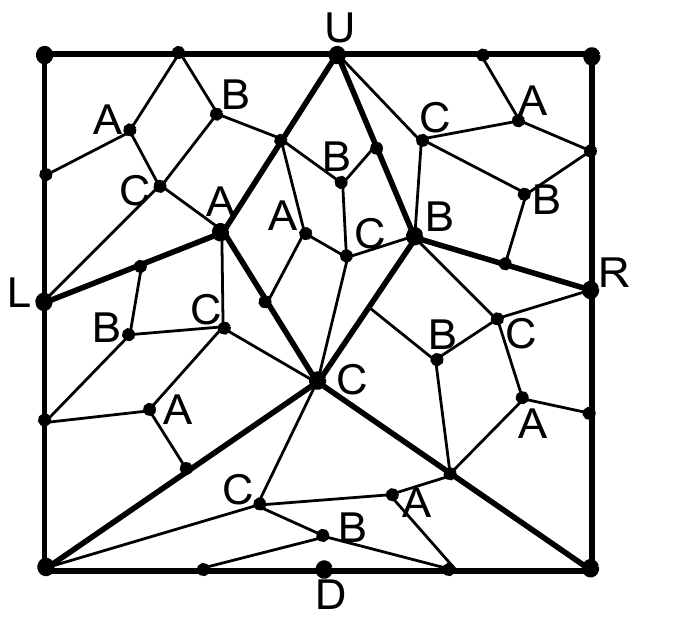}
\caption{Макроплитка третьего уровня. Отмечены типы внутренних вершин ($A$, $B$, $C$)}
\label{fig:level2}
\end{figure}

\begin{figure}[hbtp]
\centering
\includegraphics[width=0.5\textwidth]{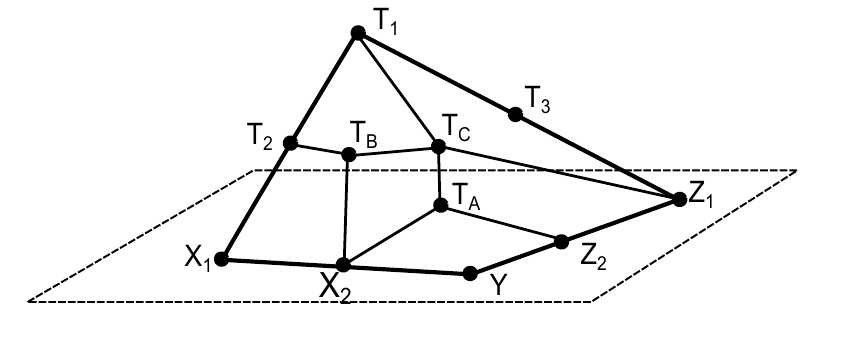}
\caption{Подклейка.}
\label{fig:pasting}
\end{figure}

Далее комплекс уровня $n$ определяется индуктивно. Сначала к комплексу $n$-уровня применяется разбиение -- каждый квадрат делится на $6$ квадратов по правилу, как на рисунке~\ref{fig:level2} потом
проводятся подклейки -- к некоторым путям длины $4$ на комплексе подклеивается $6$ квадратов. (рисунок ~\ref{fig:pasting}.

В результате образуется комплекс следующего уровня. Структура разбиений нужна чтобы обеспечить локальную преобразуемость -- любой путь с концами на границах комплекса может быть преобразован в путь только по границам локальными заменами (когда два соседних ребра некоторого квадрата меняются на два соседних ребра).

Подклейки нужны чтобы обеспечить равномерную эллиптичность комплексов: любой путь на комплексе может быть преобразован локальными заменами достаточно сильно в зависимости от его длины. Например, для естественной метрики, определенной на путях с общими концами, и некоторой константы $\alpha$, путь длины $n$ локально преобразуется в путь, отличающийся на 
$n\alpha$. 

В полугруппе, образованной кодировками путей на семействе комплексов, равномерная эллиптичность означает возможность достаточно сильно преобразововывать произвольные слова. Локальная преобразуемость же позволяет приводить слова к канонической форме.

Также нужны еще несколько свойств, позволяющих ввести раскраску всех вершин и ребер комплексов в конечное число цветов -- они будут играть роль букв в полугруппе. Также часть 
свойств посвящены геометрической структуре комплексов, помогающих приводить кодировки несуществующих на комплексе путей к нулю. 

 {\it Нулевой формой} мы называем путь по одному ребру туда-обратно.
В полугруппе слова, соответствующие {\it коротким} кодировкам перечисленных ниже типов, приравниваются к нулю:

 1) не встречающиеся на комплексах;
 
 2) пути по одному ребру туда-обратно (нулевая форма);
 
 3) принадлежащие к числу <<мертвых паттернов>> -- это пути, которые не могут быть частью достаточно длинного кратчайшего пути на комплексе.
 
 Если в процессе локальных преобразований в слове возникает короткое запрещенное подслово, мы приходим к нулю. Если путь лежит на комплексе, но является некратчайшим путем, соединяющим его концы, то можно локально преобразовать его к нулевой форме и обнулить.

Ниже приводится список основных предложений первой части работы:

\begin{lemma}[Об ограниченности роста степени вершины] \label{growth_bound}

Для каждой вершины $Z$ существует такое натуральное $N$, что
начиная с уровня макроплитки $N$, степень (число входящих ребер) вершины $Z$ не меняется, то есть она одинакова для макроплиток уровня $N$ и $N+k$ для любого натурального $k$.
\end{lemma}

{\bf Следствие.} 1. Каждая вершина заданной глубины $x$ выступает в качестве ядра подклейки для ограниченного количества вершин.

2. В каждую вершину входит ограниченное количество ребер различных уровней, включая ребра из подклеек.

Предложение нужно чтобы обосновать конечность числа возможных типов вершин. Это нужно для введения конечного алфавита полугруппы.

\begin{lemma}[О боковых и краевых вершинах] \label{side_node}

1) Каждая боковая или краевая вершина лежит на середине стороны в какой-либо макроплитке или в двух макроплитках одного уровня, лежащих в одной плоскости;

2) Если боковая или краевая вершина не находится на границе исходного комплекса первого уровня, то она лежит на одном из восьми внутренних ребер в некоторой макроплитке, либо лежит на границе подклееной макроплитки.

\end{lemma}

Предложение нужно для описания структуры кодировки типов вершин комплекса.

\begin{lemma}[О переброске пути] \label{path_flip}

Пусть $XYZT$ -- некоторая макроплитка. Рассмотрим путь $XYZ$ (состоящий из двух соседних макроребер). Тогда, если разрешается менять подпуть из двух соседних ребер любой плитки на подпуть из двух других ребер (с общими началом и концом у этих подпутей), то путь $XYZ$ может быть преобразован в $XTZ$ -- путь по другим двум макроребрам.

\end{lemma}

\begin{lemma}[О выносе пути на границу] \label{to_border}

Пусть начало и конец пути $P$, проходящего по макроплитке $T$, лежат на границе $T$. Тогда можно реализовать одну из двух возможностей:

1) $P$ может быть локально преобразован в нулевую форму.

2) $P$ может быть локально преобразован в форму $P'$ так, что $P'$ полностью лежит на границе $T$.

Кроме того, любой кратчайший путь, соединяющий противоположные углы или середины противоположных сторон макроплитки имеет длину $2^n$, и может быть локально преобразован в любой из двух путей по границе (полупериметр), с теми же концами.

\end{lemma}

Это удобные предложения для преобразований путей.

\begin{lemma}[О ``шевелении'' пути] \label{path_changing}

 Пусть $P$ -- путь, состоящий из макроребра $XY$ некоторой макроплитки $T$ уровня не менее $4$.
Обозначим середину $XY$ как $H$, середины $XH$ и $HY$ как $G_1$ и $G_2$.
 В соответствии с определением подклеек, существует подклееная макроплитка, три угла которой лежат на макроребре: это $G_1$, $H$ и $G_2$.
Тогда $P$ может быть локально преобразован в форму $P'$, состоящую из трех частей:
первая -- $XG_1$, третья -- $G_2Y$, а вторая представляет собой путь по двум соседним ребрам подклееной макроплитки, два других ребра которой являются частью $P$.
\end{lemma}

\begin{lemma}[О выделении локального участка] \label{longpath}

Пусть $n$ -- уровень макроплитки и путь $P$ лежит внутри нее.

1. Пусть оба края $P$ лежат на границе $T$. Тогда если путь $P$ имеет длину не менее $5\times 2^{n-2}$, он может быть локально преобразован в нулевую форму $P'$.

2. Пусть один край $P$ лежит в углу или на середине стороны $T$, а второй край - внутри $T$ (часть пути может проходить по границе $T$). Тогда если путь $P$ имеет длину $2^{n}$, он может быть локально преобразован в нулевую форму $P'$.

\end{lemma}

Это предложение используется для приведения к каноническому виду. Суть в том, что достаточно длинные пути на небольшом куске комплекса гарантированно являются некратчайшими.

\begin{lemma}[О непродолжаемом пути] \label{bad_path}

Рассмотрим некоторую макроплитку $T$, являющуюся подплиткой более крупной макроплитки.
Пусть путь $P$ лежит в $T$, причем начало и конец $P$ лежат в серединах противоположных сторон $T$. Тогда любые плоские пути вида $WP$ или $PW$, где длина $W$ более длины $P$, могут быть приведены к нулевой форме.
\end{lemma}

\begin{lemma}[О мертвых паттернах] \label{DeadPaterns}

Рассмотрим некоторую макроплитку $T$ и обозначим в ней внутренние вершины $A$, $B$, $C$ и боковые $U$, $R$, $D$, $L$ (аналогично обозначениям при разбиениях).
Тогда паттерны $AUB$, $ACB$, $CXD$ (где $X$ -- нижняя левая, либо нижняя правая вершина) являются мертвыми.
\end{lemma}

\begin{lemma}[О мертвых путях в нижней подплитке] \label{DeadPaths}

Рассмотрим некоторую макроплитку $T$ уровня $n$. Пусть путь $XYZ$ лежит в $T$, причем $Y$ -- левый нижний угол $T$, $XY$ лежит на внутреннем ребре, идущем из левого нижнего угла $T$ во внутреннюю вершину $C$, а $YZ$ лежит на нижней стороне $T$. Тогда для любых плоских путей $W_1$, $W_2$, длины которых более $2^{n+1}$, путь $W_1XYZW_2$ можно локально преобразовать к нулевой форме.

\end{lemma}

\begin{lemma}[О некорректных участках] \label{uncorrect_sectors}

Пусть есть некорректный участок $XYZ$ в макроплитке $T$ уровня $n$, причем $T$ -- минимальная макроплитка, содержащая $XYZ$ в качестве некорректного участка. Тогда для любых плоских путей $W_1, W_2$, длины более $2^{n+2}$, путь $W_1XYZW_2$, может быть локально преобразован к нулевой форме.
\end{lemma}

Три предложения выше разбирают суть понятия мертвого паттерна и показывают, что некоторые геометрические пути на комплексе могут быть обнулены, как не являющиеся частями длинных кратчайших путей.

\begin{lemma}[О корректности путей] \label{correct_paths}

Пусть путь $P$ представляет собой проход по двум соседним сторонам некоторой макроплитки $T$. Тогда любые локальные преобразования не могут привести $P$ к нулевой форме, а также к форме, содержащей некорректный участок.

\end{lemma}

Это предложение нужно для объяснения бесконечности полугруппы, а именно что некоторые длинные пути на комплексе не могут быть преобразованы к нулю.

\begin{lemma}[О расстоянии от края до выхода в подклейку] \label{pasting_distance}

1. Пусть вершина $X$ лежит на краю некоторой макроплитки $T$ принадлежащей комплексу $K$, вершина $Y$ принадлежит $T$, но не находится на ее границе, а из $Y$ существует выход в подклееную макроплитку уровня $n \geq 2$, не содержащую вершин на границе $T$. Тогда расстояние от $X$ до $Y$ в комплексе $K$ (длина кратчайшего пути по ребрам) не менее $2^{n-1}$.

\end{lemma}

\begin{lemma}[Об ограниченности пути, уходящего в подклееную часть] \label{pastingpath}

1. Пусть путь $P$ имеет вид $V_1AV_2BV_3CV_4$, где ребра из вершин $A$, $B$, $C$ ведут в подклееные области, а участки $V_1$, $V_2$, $V_3$, $V_4$ -- плоские. Тогда путь $P$ не может быть плоским (лежать в одной макроплитке).

2.Пусть путь $P$ начинается в вершине $X$, выход из $X$ идет в подклееную макроплитку уровня $n$. Кроме того, пусть $P$ не содержит ребер-выходов из подклееных плиток.
Тогда, если длина пути не менее $2^{n+1}$, то он может быть приведен к нулевой форме.

\end{lemma}

Эти предложения нужны для доказательства, что периодические пути с периодом $9$ гарантированно приводятся к нулю.

\medskip

\section{Содержание второй работы}

Во второй работе цикла <<Детерминированная раскраска семейства комплексов>> (\cite{coloring}) была введена раскраска вершин и полуребер данного семейства в конечное число цветов, играющих роль букв. В отдельном приложении приводится перебор, проверящий, что такая раскраска обладает свойством детерминированности, то есть если известна кодировка пути по двум ребрам, соединяющего противоположные вершины некоторого минимального квадрата (плитки) на комплексе, то определяется однозначно путь по другим двум ребрам.

Алфавит в полугруппе состоит из букв двух типов:

Реберные буквы, кодирующие всевозможные входящие и выходящие из вершин ребра на комплексах, в частности, по букве мы определяем, главное ли это ребро, входящее в вершину или не главное, и таким образом можем определять, что путь в этом месте <<вернулся>> из подклейки или подплитки, или, наоборот, <<зашел>> в подклейку или подплитку;

Вершинные буквы, кодирующие всевозможные вершины на комплексах.

В общей ситуации каждая вершина принадлежит некоторой макроплитке в базовой плоскости и также может лежать на границе каких-то подклееных макроплиток. Буквы в полугруппе кодируют 
различные варианты таких расположений. Итоговый цвет вершины состоит из комбинации следующих компонентов:

1) тип вершины;

2) уровень вершины;

3) окружение;

4) информация.

Тип вершины определяется отдельно для краевых, угловых, боковых и внутренних вершин. 

1) {\it Угловые}. (Лежащие в углах подклееных макроплиток или всего комплекса). {\it Тип угловой вершины } определим как один из четырех вариантов углов, в котором она может находиться: $\mathbb{CUL}$, $\mathbb{CUR}$, $\mathbb{CDR}$, $\mathbb{CDL}$. (Corner Up-Left и так далее.) Угловые вершины {\it принадлежат} макроплитке, где они являются углами.

2) {\it Краевые}. (Лежащие на стороне подклееной макроплитки или всего комплекса). Каждая такая вершина лежит в середине стороны некоторой макроплитки, прилегающей к краю. Тип краевой вершины определим в соответствии с тем, серединой какой стороны в этой макроплитке она является: $\mathbb{L}$, $\mathbb{R}$, $\mathbb{D}$, $\mathbb{U}$.

3) {\it Внутренние}. В этой категории определим три типа внутренних вершин: $\mathbb{A}$, $\mathbb{B}$, $\mathbb{C}$, отвечающих внутренним узлам макроплиток, эти вершины создаются {\it внутри} разбиваемой макроплитки. 

4) {\it Боковые}. (Лежащие на границе между двумя макроплитками, на внутреннем ребре).
 Типы боковых вершин будут соответствовать всем упорядоченным парам из множества $\{ \mathbb{U},\mathbb{R},\mathbb{D},\mathbb{L}\}$. А именно: $\mathbb{DR}$, $\mathbb{RD}$, $\mathbb{DL}$, и так далее. Будем считать, что в упорядоченной паре первой называется буква, соответствующая $A$-стороне внутренего ребра, на котором лежит вершина. Все боковые вершины создаются в середине стороны разбиваемой макроплитки. 
Тип боковой вершины определяет, серединой каких именно сторон она является в двух макроплитках, где она является серединой сторон. Это как раз те макроплитки, которые разбивались при создании данной вершины.

\medskip

Уровень вершины может быть 1, 2 и 3. Первый у только что созданных вершин при подклейке или разбиении (самых молодых), второй у созданных на предыдущем этапе, третий у всех остальных.

\medskip

Окружение кодирует структуру недалеко находящихся вершин. Точное определение во второй части работы (\cite{coloring}).

\medskip

{\bf Начальники.} Для каждой вершины внутренней или боковой вершин определены не более трех других вершин, считающихся начальниками данной.
Пусть $T$ -- минимальная макроплитка, такая что вершина $X$ находится внутри (не на границе) $T$.

\begin{definition}

{\it Первым начальником} вершины $X$ будем называть узел в середине верхней стороны макроплитки $T$.

Для боковых вершин, лежащих на ребрах $2$, $5$, $6$ макроплитки $T$, {\it вторым начальником} будем называть узел в правом нижнем углу макроплитки $T$.

Для боковых вершин, лежащих на ребрах $7$ и $8$, а также для для вершин типа $\mathbb{C}$ определим {\it второго} и {\it третьего начальников} как узлы в левом нижнем и правом нижнем углах $T$ соответственно.

\end{definition}

У краевых вершин начальников нет. Из угловых вершин определен начальник только для $\mathbb{CDR}$-узла, и им будет $\mathbb{CDL}$-узел в той же макроплитке.

\medskip

Информация кодирует первые три параметра (тип, уровень, окружение), но не у самой вершины $X$, а у ее начальников.

\medskip

 Указанные параметры определяются для вершины как в ее базовой плоской области, так и в подклееных областях, где она лежит, при этом не учитывается кратность, если вершина занимает одно и то же положение и получает одни и те же параметры в нескольких подклееных областях. В частности, если мы говорим,
 что в подклееном типе некоторой вершины есть $\mathbb{CUR}$, это означает, что есть подклееная макроплитка, верхний правый угол которой попадает в нашу вершину.

\medskip

Во второй части работы проводится введение описанной кодировки вершин и ребер. При этом доказывается ряд технических лемм, позволяющих вычислять параметры вершин по другим параметрам рядом расположенных вершин. Например, вычисляются параметры вершины в середине верхней стороны по параметрам вершин в углах той же макроплитки. Эти функции используются при проведении перебора: для всевозможных геометрических расположений квадрата (плитки) на комплексе проверяется, что кодировка пути по двум сторонам этой плитки однозначно определяет кодировку пути по другим двум сторонам этой плитки.

\begin{theorem}[О введении раскраски со свойством детерминированности] \label{determ}

Введенная система цветов для ребер и вершин семейства комплексов обладает свойством детерминированности: если $A$, $B$, $ C$ - вершины некоторого минимального 4-цикла $ABCD$ на комплексе, и при этом путь $ABC$ не образует мертвый патерн, то по кодировке пути $ABC$ можно однозначно определить кодировку пути $ADC$.

\end{theorem}

Эта однозначность (детерминированность кодировки) позволяет корректно ввести определяющие соотношения в полугруппе. Для каждой плитки каждые два таких пути объявляются эквивалентными, а их кодировки -- равными словами в полугруппе (если не один из них до этого не был объявлен запрещенным).

\medskip

Кроме того, ко второй части относится перебор различных случаев расположения пути по двум соседним ребрам минимального квадрата на комплексе. Перебор проверяет что кодировка парного пути, по другим двум соседним ребрам, определяется однозначно.

Задание кодировки вершин и ребер вместе с проверкой детерминированности завершает задание системы определяющих соотношений. Эта система состоит из нескольких пар эквивалентных между собой путей по двум соседним ребрам всевозможных квадратов на семействе комплексов, а также из нескольких мономиальных соотношений: 

1) слова короткой длины, не являющиеся кодировкой никакого пути на комплексе;

2) кодировки коротких некратчайших путей (пути туда и обратно по некоторому ребру);

3) кодировки мертвых паттернов -- путей, которые не могут быть продолжены до достаточно длинного пути.

Таким образом, задана конечно определенная полугруппа путей на заданном семействе комплексов. В третьей части мы доказываем, что она бесконечна (длинные пути не приводятся к нулю и не меняют длины при преобразовании), а также обладает ниль свойством (произвольное слово содержащее период 9 может быть приведено к нулю локальными преобразованиями).

\section{Приведение к канонической форме} \label{canonic}

После введения определяющих соотношений мы можем оперировать исключительно словами-кодировками путей. Учитывая, что кодировки, отвечающие коротким запрещенным путям, были обнулены, мы можем считать, что на каждом коротком подпути наш путь вкладывается в некоторый комплекс. И с ним можно проводить локальные преобразования, как с соответствующим вложенным путем. В процессе таких преобразований мы получаем либо все более крупный подпуть, вложенный в комплекс, либо приходим к запрещенному участку и нулевому слову.

Все указанные параметры (можно считать их оттенками цвета) используются при приведении слова-кодировки к каноническому виду.

Далее будет применяться терминология, введенная в первых частях работы. В частности, в цвет вершины
входит несколько параметров, таких как ее тип, окружение и информация. Цвета соответствуют буквам в полугруппе, и говоря о букве, кодирующей некоторую вершину, мы можем пользоваться ее типом, окружением и информацией в рамках алгоритма по приведению к канонической форме. 

Во второй работе были введены компоненты раскраски узлов. Фактически, мы будем использовать плоскую компоненту этой раскраски, иногда обращаясь к цвету в подклееной части. Чаще всего, нас будет интересовать, имеет ли некоторый узел в подклееной части определенный тип.

Общая схема рассуждений будет иметь вид: мы рассматриваем некоторое короткое слово, используем то, что короткие слова, не реализуемые в виде путей на комплексе, обнулены. То есть данное короткое слово можно рассматривать как кодировку некоторого пути. После этого мы применяем подходящие соотношения для преобразования пути. Далее обращаемся к другому участку. Для применения соотношений нас будет интересовать тип узлов в базовой плоскости, а также - имеет ли этот узел определенный тип в некоторой подклейке. Применение соотношений для преобразования пути всякий раз обеспечено свойством детерминированности, доказанным для всевозможных коротких путей во второй работе.


\medskip

В этой работе мы покажем, как проводить преобразования для приведения слова к канонической форме. Сначала мы рассмотрим пути, лежащие на одном ребре некоторой макроплитки, а затем разберем общий случай.

Пусть $X_i$ -- буквы, кодирующие всевозможные входящие ребра на семействе комплексов, $Y_i$ -- буквы, кодирующие узлы (их типы, окружения и информации), $Z_i$ -- буквы, кодирующие выходящие ребра. Фактически, $X_i$ и $Z_i$ выбираются из одного алфавита ребер входа-выхода.

Будем говорить, что слово $W$ имеет {\it правильную форму}, если в нем справа от любой (не последней в слове) буквы семейства $X$ обязательно стоит буква семейства $Y$, справа от любой (не последней в слове) буквы семейства $Y$ стоит буква семейства $Z$, а справа от любой (не последней в слове) буквы семейства $Z$ стоит буква семейства $X$.
Рассмотрим полугруппу $S$ с нулем $0$ и порождающими $\{X_i,Y_i,Z_i \}$.

Будем считать, что в полугруппе $S$ введены определяющие соотношения следующих категорий:

\medskip

{\bf Категория 1}: соотношения, обеспечивающие правильную форму: $X_iZ_j=0$, $X_iX_j=0$, $Y_iX_j=0$, $Z_iY_j=0$, $Y_iY_j=0$, $Z_iZ_j=0$, для всевозможных пар $(i,j)$.

\medskip

{\bf Категория 2}: $U_i=0$, где $U_i$ пробегает всевозможные слова длины не более $4$, не являющиеся кодировкой никакого пути на комплексе.

\medskip

{\bf Категория 3}: $U_i=0$, где $U_i$ пробегает всевозможные слова, являющиеся кодировкой какого-либо пути, представляющего нулевую форму (путь туда и обратно по некоторому ребру).

\medskip

{\bf Категория 4}: $U_i=0$, где $U_i$ пробегает всевозможные слова, длины не более $6$, являющиеся кодировкой какого-либо пути, представляющего путь с мертвым паттерном (одним из мертвых паттернов, перечисленных в предложении о мертвых паттернах в первой части работы).

\medskip

{\bf Категория $5$}: $U_i=V_j$, где $U_i$ и $V_j$ -- кодировки путей, участвующих в некотором локальном преобразовании. Конечное число таких пар приведено при введении определяющих сооношений во второй части работы.

\medskip

{\bf Оценка числа соотношений.} Букв в алфавите не более $7 \cdot 10^{36}$. Всевозможных слов длины не более $4$ существует не более $2{,}401 \cdot 10^{147}$. Учитывая, что в других категориях соотношений значительно меньше, можно считать, что общее число соотношений, введенных в полугруппе, не превосходит $3 \cdot 10^{147}$.

\smallskip

{\bf Замечание.} Введение соотношений можно значительно оптимизировать, сильно уменьшив их число. Но это выходит за рамки этой работы.

\medskip

Далее мы рассмотрим следствия из этих соотношений.

\subsection{Свойства путей, лежащих на одном ребре}

Пусть $W$ -- последовательность боковых узлов, принадлежащих одному внутреннему ребру некоторой макроплитки.
Поскольку мы ввели мономиальные соотношения для всех слов, не длиннее $10$, не являющихся кодировками путей на комплексе, мы можем установить, какие слова в принципе могут быть разрешены в полугруппе $S$.

\begin{definition}
Определим {\it последовательность уровней на ребре} $\mathbf{EdgeLevels}(k)$ в зависимости от натурального параметра $k$ следующим образом:

\smallskip

$\mathbf{EdgeLevels}(1)$=$1$;

$\mathbf{EdgeLevels}(2)$=$121$;

$\mathbf{EdgeLevels}(3)$=$1213121$;

Для $k>3$, $\mathbf{EdgeLevels}(k)$=$\mathbf{EdgeLevels}(k-1)$ $3$ $\mathbf{EdgeLevels}(k-1)$.

\end{definition}

Рассмотрим последовательность $\mathbf{EdgeLevels}(k)$. Заменим $1$ на $2$, $2$ на $3$, $3$ оставляем как есть. После этого в начале и конце последовательности, а также между любыми двумя членами добавим $1$. Очевидно, что при этом получится $\mathbf{EdgeLevels}(k+1)$. также можно заметить, что последовательность периодическая, с периодом $3121$ и предпериодом $121$.
Ясно, что для любого бокового ребра существует такое $k$, что последовательность $\mathbf{EdgeLevels}(k)$ отвечает значениям уровней боковых узлов на этом ребре.

\medskip

\begin{lemma} \label{edgepaths}
Пусть $W$ -- слово в $S$, в котором все буквы, соответствующие узлам, кодируют боковые узлы, лежащие на одном и том же ребре, одинаковым для всех параметром информации. Кроме того, пусть все реберные буквы одинаковы и равны $1$ или $2$ (то есть все переходы по главным ребрам). Тогда:

$1$. Уровни узлов на ребре идут в последовательности $\mathbf{EdgeLevels}(k)$ для некоторого натурального $k$.

$2$. Для каждого узла $X$, уровень которого выше $1$, существует подклееная макроплитка $T$, такая что тип $X$ относительно $T$ (тип в подклейке) равен либо $\mathbb{CUR}$ либо $\mathbb{CDL}$, причем оба таких типа (для разных макроплиток) могут быть только у одного узла в $W$, находящегося в середине внутреннего ребра макроплитки.

\end{lemma}

\begin{proof} Первое утверждение следует из того, что в любом комплексе уровни идут именно в такой последовательности. То есть, все слова длины не более $20$, не являющиеся подсловами такой периодической последовательности, присутствуют среди мономиальных соотношений. Из этого следует, что если последовательность узлов не содержит запрещенных подслов, последовательность уровней узлов соответствует указанной.

Для доказательства второго утверждения рассмотрим, как устроены последовательности узлов на внутренних ребрах построенного комплекса. После операции разбиения, узлы первого уровня становятся узлами второго. При операции подклейки, каждый узел второго уровня становится либо правым верхним, либо левым нижним углом новой подклееной макроплитки, причем в дальнейшем к этим узлам макроплитки углами не подклеиваются.
Единственным узлом, в подклееном типе которого может присутствовать и $\mathbb{CUR}$ и $\mathbb{CDL}$, является узел в середине внутреннего ребра, так как он может участвовать в нескольких подклееных макроплитках.

\end{proof}

{\bf Ранги боковых узлов}.
Пусть $W$ -- последовательность боковых узлов, принадлежащих одному внутреннему ребру некоторой макроплитки.
Будем говорить, что {\it боковой узел имеет ранг $1$} (или, соответственно, ранг $2$), если этот боковой узел -- первого (или, соответственно, второго) уровня. {\it Регулярным словом} ранга $1$ будем считать слово из одного узла ранга $1$.

Теперь индуктивно определим остальные ранги узлов. {\it Регулярным словом} ранга $n$ будем называть последовательность узлов $XYX$, где $X$ -- регулярное слово ранга $n-1$, а $Y$ -- узел ранга $n$. То есть, например, если ранги узлов идут в последовательности $121$, то это регулярное слово ранга $2$.
Будем говорить, что {\it боковой узел $X$ имеет ранг $n$}, если выполнены следующие условия:

$1.$ $X$ находится в середине подслова $\widehat{W} \subset W$, состоящего из $2^{n-1}+1$ букв;

$2.$ Слово $\widehat{W}$ представляется в форме $YUXVZ$, где $Y$, $Z$ -- боковые узлы ранга $n-1$, $U$, $V$ -- регулярные слова ранга $n-2$.

\medskip

Слово $\widehat{W}$ будем называть {\it словом-представителем} для буквы $X$ ранга $n$.

\medskip

{\it Диаграммой рангов} слова $U$ будем называть слово $\widehat{U}$ в счетном алфавите $\{ 1,2,\dots \}$, получаемое из слова $U$ заменой каждой буквы на ее ранг.
Например, узел $X$ имеет ранг $3$, если в $W$ есть подслово $Z_1Y_1XY_2Z_2$, где у $Z_1$ и $Z_2$ ранг $2$, а у $Y_1$ и $Y_2$ ранг $1$. Тогда диаграмма рангов выглядит так: $21312$.


\medskip

{\bf Замечание.}
Смысл понятия {\it слово-представитель $X$} в том что слово-представитель является половиной периметра (верхняя плюс левая сторона) подклееной макроплитки, ядром которой является $X$.

\medskip

\begin{figure}[hbtp]
\centering
\includegraphics[width=1\textwidth]{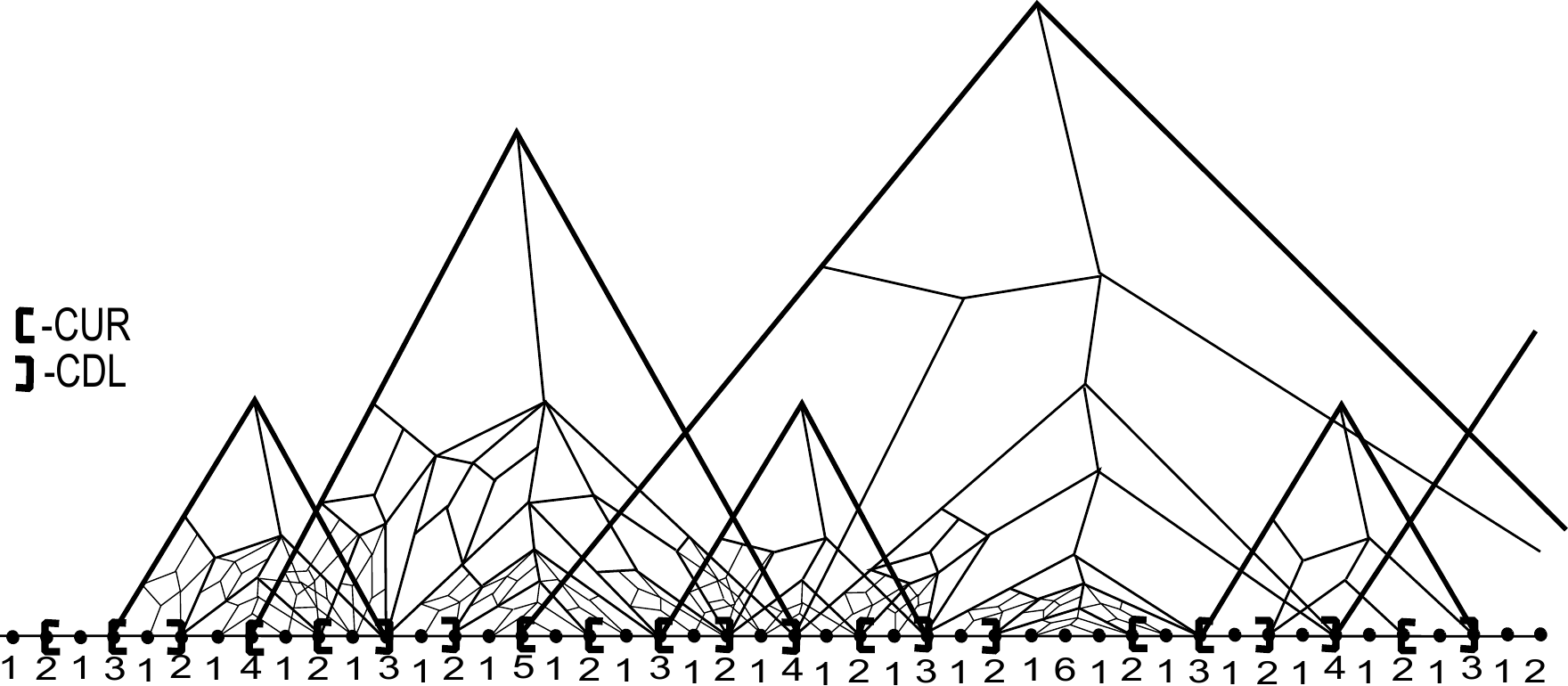}
\caption{Устройство подклеек на ребре}
\label{rankword}
\end{figure}

\medskip

\begin{lemma}  \label{skobka}
Пусть $W$ -- последовательность узлов на одном внутреннем ребре некоторой макроплитки и узел $X$ имеет ранг $n$. Рассмотрим слова длиной $2^k+1$, с центром в $X$. Его слово-представитель, имеющее длину $2^{n-1}+1$, является наименьшим среди таких слов, левый узел которых имеет в подклееном типе $\mathbb{CUR}$, а правый -- $\mathbb{CDL}$.

\end{lemma}

Это утверждение доказывается индукцией по $n$.

\medskip

\begin{lemma}  \label{mainwords}
Пусть $W$ -- слово в $S$, в котором все буквы, соответствующие узлам, кодируют боковые узлы, причем все они лежат на ребре одного и того же типа, и параметр информации у всех одинаковый.

Пусть $X$ -- буква в слове $W$, уровень которой равен $3$, причем между краями $W$ и $X$ находятся, как минимум, $4$ буквы, отвечающие кодам узлов. Тогда либо некоторое подслово $\widehat{W}\in W$, содержащее $X$, приводится к нулю с помощью определяющих соотношений в $S$, либо существует подслово в $W$, являющееся словом-представителем $X$ и оно может быть реализовано в качестве кодировки пути на некотором ребре комплекса.

\end{lemma}

\begin{proof} Будем считать, что мы читаем слова в положительном (относительно несущего ребра) направлении. Пусть уровень $X$ (то есть и ранг тоже) равен $3$. Рассмотрим подслово $W$, содержащее $5$ узловых букв, с центром в $X$.

Допустим, первая и последняя буквы содержат в подклееном типе, соответственно, $\mathbb{CDL}$ и $\mathbb{CUR}$.
Если $W$ не встречается на комплексе, оно должно присутствовать среди мономиальных определяющих соотношений. Иначе оно автоматически реализуется в качестве кодировки пути.
 В этом случае $W$ -- кодировка пути по верхней и левой стороне подклееной макроплитки уровня $2$, и $W$ -- слово-представитель $X$.

Пусть, теперь, первая и последняя буквы содержат в подклееном типе, соответственно, $\mathbb{CUR}$ и $\mathbb{CDL}$, то есть в другом порядке. (Иные варианты не встречаются на комплексе в качестве кодировок и потому приводятся к нулю).

\begin{figure}[hbtp]
\centering
\includegraphics[width=1\textwidth]{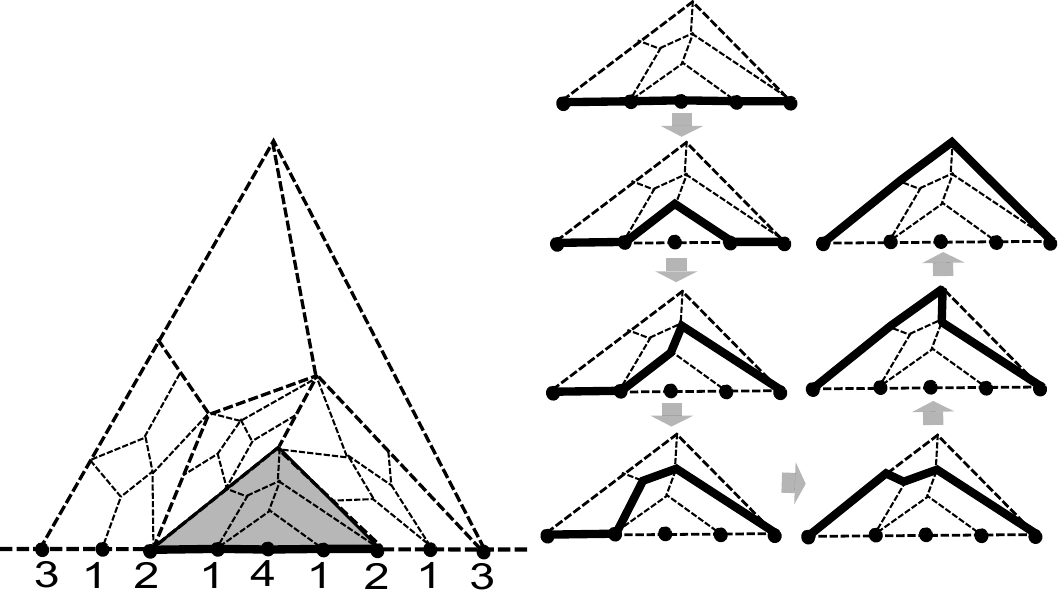}
\caption{Преобразование слова-представителя}
\label{enlarge}
\end{figure}

В этом случае ранг центральной буквы -- $4$ или выше. Рассмотрим слово $V_1$ из девяти букв с центром в $X$. Диаграмма рангов $V_1$ имеет вид $3121x1213$. Слово $V_1$ встречается на комплексе, иначе все можно сразу привести к нулю. Рассмотрим типы (в подклееной части) для четвертой и шестой букв этого слова. Это должны быть разные типы, а учитывая положительное прочтение слова, четвертая должна иметь тип $\mathbb{U}$, а шестая $\mathbb{L}$ (речь о типах в подклееной области). То есть, к подслову из четвертой, пятой и шестой букв может быть применено соотношение для локального преобразования из числа введенных в случае $P2$. Затем можно применить другие локальные преобразования $7$ из числа введенных в случае $P2$, для соответствующих трехбуквенных слов. В результате мы преобразуем слово $V_1$ в кодировку пути, изображенного на рисунке~\ref{enlarge}.

\medskip

\begin{figure}[hbtp]
\centering
\includegraphics[width=0.9\textwidth]{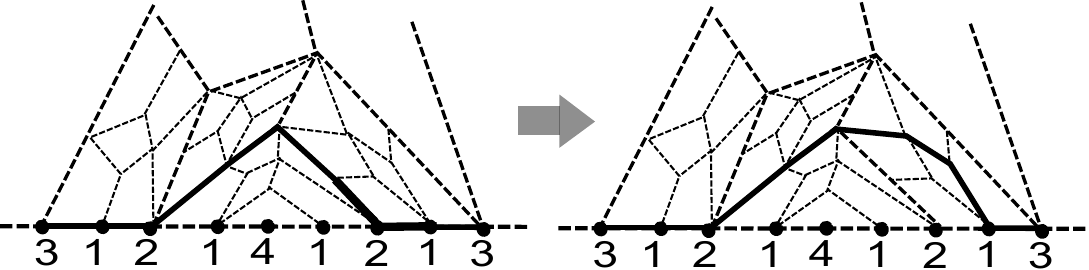}
\caption{Преобразование слова-представителя}
\label{enlarge2}
\end{figure}

\medskip

Рассмотрим трехбуквенное подслово, отмеченное в левой части на рисунке~\ref{enlarge2}. Заметим, что либо путь с таким кодом не может встретиться на комплексе (и тогда оно равно нулю), либо к нему можно применить соотношение $P3-8$. Действительно, узел $X$ обязан иметь уровень $1$, значит, базовый тип его $\mathbb{UL}$ или $\mathbb{LU}$, а тип в подклейке -- $\mathbb{U}$ (если будет $\mathbb{L}$, то такого участка из трех узлов не может встретиться на комплексе). Затем, аналогично применяем соотношения $B3-6$, $P3-1$, и получаем кодировку пути, изображенного в правой части рисунка~\ref{enlarge2}.

Далее возможно два варианта. Последняя буква слова $V_1$ может содержать в подклееном типе либо $\mathbb{CDL}$ либо $\mathbb{CUR}$.

В первом случае мы можем применить соотношение $P3-10$, далее $B3-5$, $B3-2$. Применяя далее подходящие соотношения, можно получить из слова кодировку пути, изображенного в левой части на рисунке~\ref{enlarge3}. Согласно предложению~\ref{skobka}, в этом случае рассмотренное девятибуквенное слово и есть слово-представитель для $X$. Применив далее несколько соотношений, мы можем преобразовать это слово в кодировку пути, указанного в правой части рисунка~\ref{enlarge3}.

\begin{figure}[hbtp]
\centering
\includegraphics[width=0.9\textwidth]{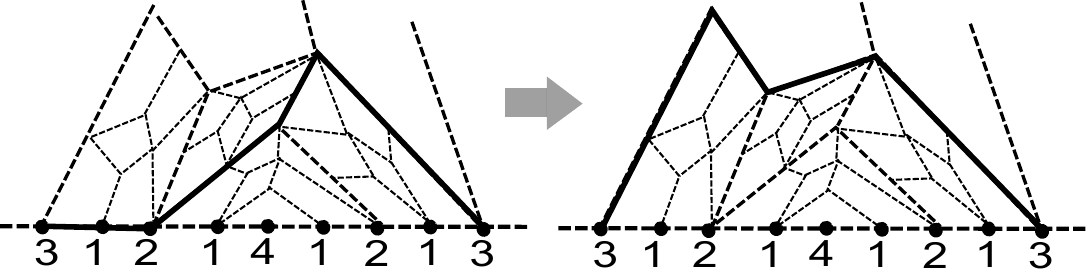}
\caption{Преобразование слова-представителя}
\label{enlarge3}
\end{figure}

\smallskip

Во втором случае мы продолжаем применять соотношения, пока не получим из нашего девятибуквенного слова кодировку указанного на рисунке~\ref{enlarge4} пути.

\begin{figure}[hbtp]
\centering
\includegraphics[width=0.7\textwidth]{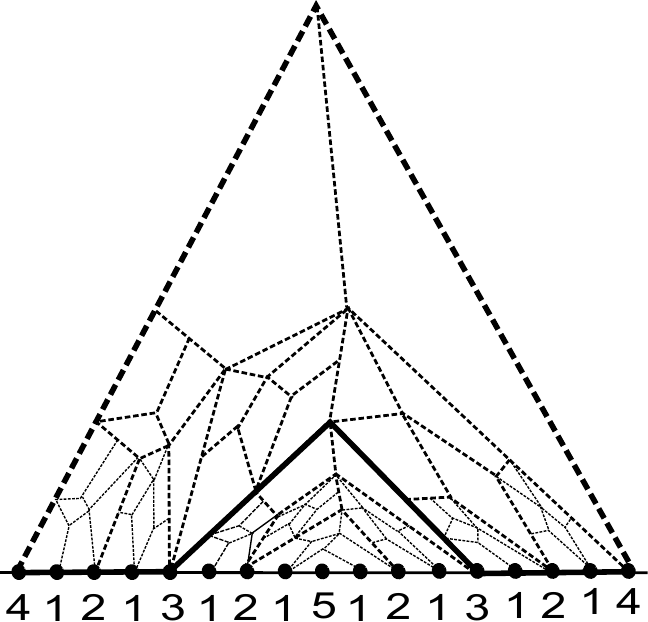}
\caption{Преобразование слова-представителя}
\label{enlarge4}
\end{figure}

Теперь рассмотрим $17$-буквенное слово $V_2$ с центром в $X$ (рисунок~\ref{enlarge4}). Мы опять можем брать трехбуквенные слова на границе $17$-буквенного слова так, что либо эти трехбуквенные слова будут встречаться в мономиальных соотношениях, либо в соотношениях описанных в главах про локальные преобразования путей. После проведения необходимых преобразований, мы опять получим альтернативу: последняя буква слова $V_2$ может содержать в подклееном типе либо $\mathbb{CDL}$ либо $\mathbb{CUR}$. Аналогично $9$-буквенному слову, мы либо получаем кодировку пути, полученного из слова-представителя, либо возможность продолжить процесс и перейти к рассмотрению $33$ буквенного слова.

Продолжая процесс преобразования, мы дойдем до того момента, когда сможем окончательно выписать кодировку пути, полученного из слово-представителя, либо получим в какой-то момент нулевое слово.

\end{proof}

\medskip

\begin{lemma}  \label{mainwords2}

Пусть $W$ -- слово в $S$, удовлетворяющее следующим условиям:

$1$. $W$ состоит из трех частей $W=V_1XV_2$, где

$X$- узловая буква, соответствующая внутреннему узлу ($\mathbb{A}$, $\mathbb{B}$ или $\mathbb{C}$), либо $X$ -- любая узловая буква, но ребро входа в $X$ или ребро выхода из $X$ -- неглавное;

$V_1$,$V_2$ -- слова, все узловые буквы которых обладают одной и той же информацией (своей для $V_1$ и для $V_2$), а все реберные буквы (кроме, возможно, входящего в $X$ и выходящего из $X$) соответствуют главным ребрам, причем тип ребра постоянен для $V_1$ и для $V_2$.

Тогда либо некоторое подслово $\widehat{W}\in W$, содержащее $X$, приводится к нулю с помощью определяющих соотношений в $S$, либо существует подслово в $W$ являющееся словом-представителем $X$ и оно может быть реализовано в качестве кодировки пути на некотором ребре комплекса.

\end{lemma}

Доказательство полностью аналогично лемме~\ref{mainwords}. В данном случае в процессе доказательства используются соотношения, связанные с подклееной макроплиткой с ядром в $X$.

\medskip

\begin{lemma} \label{checkedges}

Пусть $W$ -- слово в $S$ и каждая узловая буква $W$ кодирует боковой узел, причем информация всех таких узлов совпадает, а каждая реберная буква соответствует главному ребру (то есть для всех $1$ или для всех $2$). Тогда либо $W$ можно привести к нулю с помощью определяющих соотношений в $S$, либо $W$ представляется в виде кодировки некоторого пути на комплексе.

\end{lemma}

\begin{proof} Заметим, что если уровень всех узловых букв не выше $2$, то $W$ не может быть содержать более $3$ узловых букв, в этом случае $W$ имеет малую длину и если в нем не содержится запрещенных слов, значит оно реализуется в качестве кодировки пути.

Выберем некоторую букву $X$, уровня $3$. Пусть ее ранг равен $k$. Применим к $X$ предложение~\ref{mainwords}. Получим $V(X)$ -- слово-представитель $X$. Обозначим крайнюю справа букву $V(X)$ как $X_1$. Ее ранг будет меньше на $1$, чем ранг $X$. Применим к $X_1$ предложение~\ref{mainwords}, опять возьмем в $V(X_1)$ крайнюю справа букву. Продолжим так применять предложение~\ref{mainwords}, пока ранг и уровень всех буквы $X_i$,  которой равны $3$. Если в какой-то момент слово $W$ закончилось, мы можем рассмотреть аналогичный процесс, взяв на каждом шаге крайние слева буквы слов-представителей.


\medskip

Будем пока считать, что слово не заканчивается. Тогда мы в итоге получим букву $X_i$, ранг и уровень которой равны $2$. Следующая за $X_i$ узловая буква должна иметь уровень $1$, иначе возникнет запрещенное слово. Рассмотрим еще одну следующую узловую букву $Y$. Ее уровень равен $3$, применим к ней предложение~\ref{mainwords}. Можно заметить, что левый конец полученного слова-представителя $V(Y)$ -- это буква $X$. То есть, ранг $Y$ на один выше, чем ранг $X$.

Таким образом, с помощью указанного процесса, мы можем по заданной узловой букве либо найти букву с большим рангом, либо дойти до конца слова (в обе стороны). Во втором случае указанная буква будет обладать наивысшим рангом из всех узловых букв.

Пусть буква $X$ -- обладает максимальным рангом из всех узловых букв $W$. Теперь применим лемму~\ref{mainwords} к $X$ и рассмотрим крайние буквы $X_L$ и $X_R$ полученного слова-представителя. Заметим, что их ранг на $1$ меньше, чем ранг $X$. Для каждой из этих букв зафиксируем пару соседних с ней реберных букв так, что одно ребро уходит в подклейку (соответствующую слову-представителю $X$), а другая -- соответствует главному ребру. Применим к этим буквам с такими парами ребер предложение~\ref{mainwords2}, получая расширение подслова, представимого в качестве кодировки некоторого пути. Продолжим далее аналогичным образом, пока крайние буквы не будут иметь ранг $2$. Следующие две буквы с обоих концов должны иметь ранг $1$, и полученное таким образом слово представляется в виде кодировки пути на комплексе. Заметим, что если после этого слева или справа имеется какая-либо буква, она должна иметь больший ранг, чем $X$ (в этом можно убедиться, применив для этой буквы лемму~\ref{mainwords}). Таким образом, все слово может быть представлено в виде кодировки пути на комплексе.

Пусть у нас получился путь $P$ на комплексе. Теперь можно проделать над ним все локальные преобразования в обратном порядке, возвращая слову форму, где все узловые буквы принадлежат одному ребру и все переходы -- главные. Путь $P$ преобразуется в эквивалентный путь, идущий по одному внутреннему ребру некоторой макроплитки.

\end{proof}

\begin{lemma} \label{checkedges2}

Пусть $W$ -- слово в $S$, удовлетворяющее следующим условиям:

$W$ состоит из трех частей $W=V_1XV_2$, где

$X$ -- узловая буква, соответствующая внутреннему узлу ($\mathbb{A}$, $\mathbb{B}$ или $\mathbb{C}$), либо $X$ -- любая узловая буква, но ребро входа в $X$ или ребро выхода из $X$ -- неглавное;

$V_1$, $V_2$ -- слова, все узловые буквы которых обладают одной и той же информацией (своей для $V_1$ и для $V_2$), а все реберные буквы (кроме, возможно, входящего в $X$ и выходящего из $X$) соответствуют главным ребрам, причем тип ребра постоянен для $V_1$ и для $V_2$.

Тогда либо $W$ можно привести к нулю с помощью определяющих соотношений в $S$, либо $W$ представляется в виде кодировки некоторого пути на комплексе.

\end{lemma}

Доказательство аналогично лемме~\ref{checkedges}, учитывая что вершина $X$ тоже обязана быть ядром подклееной макроплитки со сторонами на $V_1$ и $V_2$.

\medskip

\begin{lemma}[О пути внутри макроплитки] \label{innerpath}

Пусть $W$-- слово в $S$, отвечающее следующим условиям:

{\bf $1$.} В слове $W$ отсутствуют буквы, кодирующие входы-выходы в подклейку;

{\bf $2$.} Первая и последняя буквы слова $W$ -- узловые;

{\bf $3$.} Вторая буква слова $W$ кодирует неглавное выходящее ребро;

{\bf $4$.} Предпоследняя буква слова $W$ кодирует неглавное входящее ребро;

{\bf $5$.} Кроме второй и предпоследней букв в слове $W$ ни одна буква не кодирует неглавных входящих или выходящих ребер (то есть все остальные буквы, кодирующие ребра -- отвечают главным ребрам).

Тогда $W$ либо может быть приведено к нулю, либо $W$ является кодировкой некоторого пути, соединяющего два узла на периметре некоторой макроплитки.

\end{lemma}

\begin{proof} Можно считать, что слово $W$ имеет правильную форму, то есть буквы, кодирующие входящие ребра, узлы и выходящие ребра, встречаются в слове в правильном порядке. В противном случае слово сразу можно привести к нулю, c помощью соотношений категории $1$. Также, будем считать что всякое подслово длины $4$ в нашем слове является регулярным, иначе, можно применить мономиальное соотношение из категории~$2$.

Учитывая правильность формы, буквы $1$, $4$, $7$, $\dots$, $3k+1$ слова $W$ кодируют узлы. Буквы $2$, $5$, $8$, $\dots$, $3k+2$ -- кодируют выходящие ребра, а буквы $3$, $6$, $9$, $\dots$, $3k$ -- кодируют входящие ребра.

Неглавное ребро входа означает, что путь входит внутрь некоторой макроплитки.
Заметим, всего существует шесть возможных ``входов'' внутрь макроплитки: по ребрам $1$ и $2$ (из середины верхней стороны), по ребру $3$ из середины левой стороны, по ребру $6$ из середины правой стороны, по ребру $7$ из левого нижнего угла и по ребру $8$ из правого нижнего угла.

\medskip

По второй букве $W$, кодирующей вход в макроплитку, можно установить, с каким из указанных шести входов мы имеем дело. Допустим, сначала, что это вход по ребру $3$ из середины левой стороны, а выход (предпоследняя буква) -- по ребру $1$ в середину верхней стороны.

Допустим, четвертая буква $W$ соответствует черному внутреннему узлу ($\mathbb{A}$, $\mathbb{B}$ или $\mathbb{C}$). Рассмотрим седьмую букву, соответствующую третьему по очередности узлу вдоль нашего пути. Если эта буква является последней в $W$, то $W$ представляет собой путь из середины левой стороны через узел типа $\mathbb{A}$ в середину верхней стороны. Все слова не соответствующие этому, не кодируют никакого пути на комплексе и поэтому содержатся среди мономиальных соотношений категории $2$. Значит, седьмая буква не является последней. Заметим, что если она соответствует боковому узлу, то подслово $W$ из первых семи букв будет кодировать несуществующий путь. В случае, когда седьмая буква соответствует внутреннему узлу, можно рассмотреть возможные типы десятой и тринадцатой букв. Можно заметить, что если путь с такой кодировкой существует, в нем содержится участок, проходимый туда и обратно (нулевая форма) и в этом случае, в нашем слове будет подслово из категории $3$ мономиальных соотношений.

\medskip

Теперь пусть четвертая буква $W$ кодирует боковой узел. Пусть в слове $W$ все буквы кодируют боковые узлы. Участок пути из одного бокового узла в другой по главному ребру не меняет тип этого ребра. В нашем случае это ребро $3$, причем с левой стороны -- сторона $\mathbf{3A}$, а с правой -- $\mathbf{3B}$. Но при выходе (когда мы приходим в середину верхней стороны) ребро должно быть типа $1$, то есть найдется локальный участок, где тип ребра меняется, значит можно применить соотношение из категории $2$.

Итак, после цепочки боковых узлов, нам должен встретиться внутренний узел. Поскольку у нас ребро $3$, это может быть только узел типа $\mathbb{A}$, для других не будет реализующих путей на комплексе. После узла $\mathbb{A}$, мы можем продолжать путь по ребру $1$, либо по ребру~$4$. (Если путь далее продолжается по ребру $3$, то образуется участок с нулевой формой и опять же можно применить мономиальное соотношение.)

В случае ребра $4$, мы продолжим изучение последующих букв, кодирующих узлы. Аналогично можно установить, что после нескольких боковых узлов, мы должны встретить узел типа $\mathbb{C}$, иначе можно будет применить мономиальное соотношение. После узла $\mathbb{C}$, дальнейший путь может продолжиться по одному из ребер $7$, $8$ или $5$. В первых двух случаях, если после цепочки боковых, мы встретим внутренний узел, можно будет применить мономиальное соотношение, так как такого участка пути не может существовать на комплексе.

Если же дальнейший путь идет по ребру $5$, опять рассмотрим цепочку боковых узлов до первого внутреннего узла (опять же, боковые узлы не могут поменять тип ребра, а в конце ребро должно быть первого типа). Это может быть только узел типа $\mathbb{B}$ и дальнейший путь может идти только по ребрам $2$ или $6$. В обоих случаях, для любого встреченного узла внутреннего типа можно будет применить мономиальное соотношение.

\medskip

Таким образом, наше слово должно содержать цепочку боковых узлов на ребре $1$, потом $\mathbb{A}$-узел, потом цепочка боковых узлов типа $3$. Если указанные цепочки содержат не более $10$ букв, кодирующих узлы, то можно напрямую либо найти все наше слово среди запрещенных, а если его там нет, значит слово реализуется в качестве кодировки. Пусть указанные цепочки достаточно длинные. Временно выбросим из слова три первые и три последние буквы. Останется слово, к которому мы можем применить лемму~\ref{checkedges2}, выбрав в качестве $X$ наш $\mathbb{A}$-узел.

То есть, урезанное слово без первых трех и последних трех букв может быть представлено в качестве кодировки пути на комплексе. Пусть это будет путь внутри некоторой макроплитки $T$, с заданными окружениями тех узлов на границе $T$, которые являются начальниками узлов из $W$.

При этом первые три и последние три буквы кодируют вход и выход в макроплитку, плюс информацию двух узлов на границе. Если вход или выход не соответствует ребрам, по которым далее должен идти путь, мы сразу можем применить мономиальное соотношение. Значит, окружение макроплитки и окружение и информация двух входных узлов могут сочетаться в одном объекте на комплексе.
Тогда получаем, что наше слово реализуется в качестве кодировки некоторого пути на комплексе.

\end{proof}

{\bf Замечание $1$}. Если вход и выход осуществляются не по ребрам $3$ и $1$, можно провести полностью аналогичный побуквенный анализ слова, и найти подслово категорий $2$ или $3$, либо установить, что слово есть кодировка некоторого пути в макроплитке. Таким образом, предложение можно доказать для входов внутрь макроплитки по ребрам любых внутренних типов.

\medskip
{\bf Замечание $2$}. Неглавные ребра в начале и конце слова могут также быть ребрами в подклейку, рассуждения при этом не меняются.

\begin{lemma}[О сокращении подпути по подклееной макроплитке] \label{deletepath}

Пусть $W$ -- слово в $S$, причем первая и последняя буквы кодируют некоторые узлы (то есть это некоторые из букв $Y_i$), а вторая и предпоследняя -- кодируют, соответственно входящее ребро в подклейку и выходящее ребро из подклейки, причем кроме этих двух, больше входов и выходов в подклееные области в слове $W$ не встречаются. Тогда $W$ либо может быть приведено к нулю, либо к форме, содержащей менее двух входов или выходов в подклееные области.

\end{lemma}

\begin{proof} Можно считать, что $W$ имеет правильную форму, то есть буквы, кодирующие входящие ребра, узлы и выходящие ребра, встречаются в слове в правильном порядке.

{\it Cкобочной структурой} слова будем называть слово в алфавите из двух букв:$\{ \mathbf{[}$, $\mathbf{]}\}$  открывающей и закрывающей скобок, составленное следующим образом: -- каждая буква, кодирующая выходящее неглавное ребро, или ребро, ведущее в подклееную макроплитку, заменяется на ``$[$'', а каждая буква, кодирующая входящее неглавное ребро, или ребро, ведущее из подклееной макроплитки, заменяется на ``$]$''. Остальные буквы заменяются пустыми буквами.

\medskip

Рассмотрим скобочную структуру $W$. Она начинается открывающей скобкой и заканчивается закрывающей. Очевидно, что существует подслово в $W$, которому в скобочной структуре соответствует подслово ``$[$ $]$'' (из открывающей и закрывающей скобок). Это подслово $\widehat{W}$ удовлетворяет условиям леммы~\ref{innerpath}, то есть, в некоторой макроплитке есть путь, соединяющий точки периметра, с таким кодом. Но тогда есть цепочка локальных преобразований, приводящая этот путь к эквивалентному, но идущему по периметру той же макроплитки. Поскольку локальным преобразованиям отвечают соотношения из категории $5$, подслово $\widehat{W}$ можно привести к форме, скобочная структура которой содержит не более одной скобки (а не две). Отметим также, что начальная и конечная скобки не исчезают при такой операции.

Будем проводить такие операции, пока во всем слове $W$ не останется только открывающей и закрывающей скобок. Теперь наше слово соответствует пути, все ребра которого -- главные, кроме первого и последнего. Данная ситуация аналогична условиям предложения~\ref{innerpath}. Таким образом, если наше слово ненулевое, оно представляет некоторый путь и мы можем провести цепочку локальных преобразований, чтобы этот путь шел по периметру макроплитки. Такой путь не может содержать более одного ребра в подклейки. Таким образом, в результате преобразований, в слове понизилось количество ребер в подклейки.

\end{proof}

\begin{lemma}[О двух выходах в подклейку подряд] \label{twoexits}

1.(Плоский случай). Пусть $W$-- слово в полугруппе $S$, причем в $W$ не входит букв, кодирующих входы и выходы в подклейки, а также кодирующих выходы из подплиток.
Пусть также первая буква является узловой, а вторая кодирует ребро в подплитку.
Тогда либо $W$ может быть приведено к нулю в $S$, либо $W$ является кодировкой некоторого пути на комплексе.

\smallskip
2. Пусть $W$-- слово в полугруппе $S$, причем в $W$ не входит букв, кодирующих выходы в подклейки, а также кодирующих выходы из подплиток. Пусть также первая буква является узловой, а вторая кодирует ребро в подклейку.
Тогда либо $W$ может быть приведено к нулю в $S$, либо $W$ является кодировкой некоторого пути на комплексе.

\end{lemma}

\begin{proof}
Узловые буквы, соответствующие узлам $\mathbb{A}$, $\mathbb{B}$, $\mathbb{C}$, $\mathbb{R}$, $\mathbb{D}$, $\mathbb{CDR}$ будем называть {\it особыми}, как и соответствующие узлы.

1. Выделим в слове $W$ буквы, кодирующие переход в подплитку. Пусть буквы $Z_i$ -- это узловые буквы, входящие в слово перед ними. Тогда слово представляется в виде $Z_1V_1\dots Z_kV_k$, где в словах $V_i$ первая буква кодирует вход в подплитку, а все остальные реберные буквы кодируют главные ребра. Зафиксируем $i$ и рассмотрим участок $Z_iV_i$. Выделим в нем особые буквы (в данном случае, они соответствуют узлам $\mathbb{A}$, $\mathbb{B}$, $\mathbb{C}$). Представим $Z_iV_i$ в виде $Z_iU_1S_1\dots S_{l-1}U_lZ$, где $U_j$ -- слова без особых букв, $S_j$ -- особые буквы. Применяя лемму~\ref{checkedges2} к участку $Z_iU_1S_1U_2$, мы получаем путь по двум внутренним ребрам некоторой макроплитки, начинающийся на ее границе. Будем говорить, что макроплитка, в которой путь с заданной кодировкой может проходить, {\it удовлетворяет условию}. Заметим, что на комплексе может существовать несколько макроплиток, удовлетворяющих условию, например, если ни у одного узла на пути нет третьего начальника, то третьим начальником у макроплитки могут быть разные узлы.

Применяя эту же лемму~\ref{checkedges2} к каждому участку $U_jS_jU_{j+1}$, мы также получаем набор макроплиток, где мог бы располагаться этот участок из двух ребер. Докажем, что множество макроплиток, удовлетворяющих всем условиям, не является пустым. Фактически, достаточно проверить возможность выбора трех начальников макроплитки, обеспечивающих все пути согласованными кодировками. Заметим, что все участки состоят из двух внутренних ребер. Но если часть участков требует одного начальника, а часть другого, то можно выбрать один участок, где эти требования не согласуются, тогда данный участок нельзя вложить в макроплитку, что является противоречием. Таким образом, существует некоторое непустое множество троек начальников, удовлетворяющих всем условиям.

Итак, весь участок $Z_iV_i$ можно представить как кодировку пути в некоторой макроплитке $T_i$. Также, участок $Z_{i+1}V_{i+1}$ также кодирует путь в некоторой макроплитке $T_{i+1}$, причем $T_{i+1}$ должна занимать определенное кодом $Z_iV_i$ положение одной из подплиток $T_i$. Допустим, опять, множества возможных положений для $T_{i+1}$ определяемых кодами $Z_iV_i$ и $Z_{i+1}V_{i+1}$ не пересекаются. Но получаем противоречие с применением леммы~\ref{checkedges2} для участка, одна часть которого лежит в $V_i$, а другая -- в $V_{i+1}$.

Таким образом, существует некоторая макроплитка, содержащая весь путь с кодировкой $W$.

\medskip

2. Выделим в слове $W$ буквы, кодирующие переход в подклейку. Пусть буквы $Z_i$ -- это узловые буквы, входящие в слово перед ними. Тогда слово представляется в виде $Z_1V_1\dots Z_kV_k$, где в словах $V_i$ первая буква кодирует вход в подклейку, а все остальные реберные буквы кодируют плоские ребра. Для участков $Z_iV_i$ можно применить пункт 1) настоящей леммы. Аналогично рассуждению выше, существуют макроплитки для участков $Z_iV_i$ и $Z_{i+1}V_{i+1}$ вместе. Иначе опять же можно применить лемму~\ref{checkedges2} для участка, одна часть которого в одном пути, а другая в другом. В итоге получаем макроплитку, где содержится общий путь для $W$.




\end{proof}

\begin{lemma}[О безквадратности пути по одному ребру]  \label{nosquares}

Пусть $W$-- путь на комплексе, целиком проходящий внутри некоторого ребра (внутреннего или граничного) и не содержащий узлов с типами $\mathbb{DL}$,$\mathbb{LD}$, $\mathbb{DR}$, $\mathbb{RD}$, $\mathbb{RU}$, $\mathbb{UR}$, $\mathbb{R}$, $\mathbb{D}$.
Тогда соответствующее слово в полугруппе $S$ не содержит двух одинаковых подслов, идущих подряд.
\end{lemma}

\begin{proof} Рассмотрим подстановочную систему (или DOLL-систему) в алфавите $\{ U_1,L_1, U, L \}$, определяемую следующим образом: $U_1\rightarrow U_1UL_1$, $L_1\rightarrow U_1LL_1$, $U\rightarrow U$, $L\rightarrow L$.

Такая система, в частности, генерирует следующую последовательность слов:

$$U_1\rightarrow U_1UL_1 \rightarrow  U_1UL_1 U U_1LL_1 \rightarrow U_1UL_1 U U_1LL_1 U U_1UL_1 L U_1LL_1 \rightarrow \dots $$

Пусть $X$ -- некоторое слово. Обозначим описанную подстановочную замену как $f(X)$.

\medskip

Пусть сначала путь отвечающий $W$ проходит по внутреннему ребру. Рассмотрим кодировку нашего слова $W$ в полугруппе $S$. Проведем над ним следующие операции:

\medskip

{\it каждую букву, кодирующую $\mathbb{UL}$-узел, уровня выше первого, заменим на букву $U$;}

{\it каждую букву, кодирующую $\mathbb{LU}$-узел, уровня выше первого, заменим на букву $L$;}

{\it каждую букву, кодирующую $\mathbb{UL}$-узел первого уровня, заменим на букву $U_1$;}

{\it каждую букву, кодирующую $\mathbb{LU}$-узел первого уровня, заменим на букву $L_1$;}

{\it все буквы, кодирующие не узловые буквы, удалим из слова.}

\medskip

Получившееся слово в алфавите $\{ U_1,L_1, U, L \}$ будем называть $W'$. Заметим, что если $W$ содержит два одинаковых подслова, идущих подряд, то таким же свойством будет обладать и $W'$.

Покажем, что $W'$ является подсловом $f^n(U_1)$ или $f^n(L_1)$ для некоторого $n$. Действительно, выберем в нашем слове $W$ букву, $X$ кодирующую узел с самым большим рангом $n$ (узел, появившийся ранее всех при разбиениях). Такая буква только одна, так как если бы их было хотя бы две, то согласно правилам разбиения, между ними нашелся бы узел с более высоким рангом.

Допустим, ради определенности, что узел $X$ имеет тип $\mathbb{UL}$. Можно заметить, что последовательность узловых букв в слове $W'$ является подсловом в слове $f^{n-1}(U_1)$, это легко можно установить индукцией по $n$, подстановочная последовательность $f$ полностью соответствует образованию узловых вершин при разбиениях.

\medskip

Таким образом, для доказательства предложения, достаточно показать, что $f^{n}(U_1)$ не содержит двух одинаковых подслов подряд для любого $n$.

Допустим, это не так. Возьмем наименьшее $n$, при котором такие подслова появляются.
Пусть $f^{n}(U_1)$ содержит подслово $QQ$. Разметим наше слово следующим образом. Выделим в $f^{n}(U_1)$ обычными (не квадратными) скобками участки, соответствующие подсловам $U_1UL_1$, $U_1LL_1$, $U$, $V$, образовавшимся при последнем переходе
$f^{n-1}(U_1) \rightarrow f^{n}(U_1)$. Сами скобки не принадлежат алфавиту, а являются просто инструментом разметки. Если при этом в $Q$ вошло целое количество участков, начинающихся и заканчивающихся скобкой,
то в слове $f^{n-1}(U_1)$ также будет два одинаковых подслова $f^{-1}(Q)$ идущих подряд, что дает противоречие с выбором $n$.

\medskip

Значит есть подслово $(U_1UL_1)$ или $(U_1LL_1)$, часть которого лежит вне $Q$, а часть внутри. То есть, $Q$ начинается с $UL_1)$ или $LL_1)$, либо с $U_1)$ или $L_1)$.

Рассмотрим первый случай. $Q$ должно заканчиваться на $(U_1$, то есть $Q$ представляется в виде $UL_1)Y(U_1$, где $Y$ -- слово содержащее целое количество участков в скобках и к нему можно применить $f^{-1}$. Тогда $f^{n}(U_1)$ представляется в виде
$$W_1 (U_1 \bf{UL_1)Y(U_1UL_1)Y(U_1} UL_1) W_2,$$

где $W_1$, $W_2$ -- некоторые подслова, содержащие целое число участков в скобках.

Можно заметить, что при применении $f^{-1}$ к подслову $(U_1 \bf{UL_1)Y(U_1UL_1)Y(U_1} UL_1)$ там также образуется два одинаковых подслова идущих подряд. То есть такие слова встретятся в $f^{n-1}(U_1)$, что невозможно.

Остальные случаи: $Q$ начинается с $LL_1)$; $U_1)$; $L_1)$ рассматриваются полностью аналогично.

\end{proof}

\medskip

\begin{lemma} \label{powerpath}
Пусть $W$ -- произвольное слово в $S$. Тогда слово $W^9$ приводится к нулю.

\end{lemma}

\begin{proof} Отметим в слове $W$ открывающими скобками -- ребра, кодирующие входы в подклейки и входы в подплитки (неглавные ребра), а закрывающими -- выходы из подклеек и подплиток. Если подслово начинается с открывающей скобки, заканчивается закрывающей, и не содержит внутри скобок, будем применять предложение~\ref{deletepath}. Таким образом, в какой-то момент мы преобразуем наше слово $W$ в форму $\widehat{W}=W_1W_2$, где $W_1$ содержит только закрывающие скобки, а $W_2$ только открывающие.

Заметим, что $W^{n+1}=W_1 (W_2W_1)^n W_2$.
Рассмотрим слово $W_2W_1$. Применяя к нему такой же процесс (предложение~\ref{deletepath} к каждой паре скобок), мы в итоге получим слово $V$, либо не содержащее скобок вообще, либо содержащее только один их вид.

Пусть, например, $V$ содержит только открывающие скобки.

Заметим, что $V^n$ содержит не менее $8$ открывающих скобок, то есть содержит подслово $U=Z_1x_1V_1 Z_2x_2V_2\dots Z_kx_kV_k$, где $x_i$ -- ребра, кодирующие выходы в подклейки, участки $V_i$ не содержат таких букв, $Z_i$ -- узловые буквы, $k\ge 8$.

Учитывая лемму~\ref{twoexits} $U$ может быть представлено как кодировка некоторого пути по комплексу. Пусть $s$ -- уровень подклееной макроплитки $T$, куда уходит ребро $x_3$. По первой части леммы~\ref{pastingpath}, кусок пути от $Z_1$ до $Z_3$ не может лежать в одной макроплитке. То есть $Z_1$ не лежит в $T$. Тогда длина пути от $Z_1$ до $Z_3$ не менее $2^{s-1}$ по предложению о пути по подклееной части подплитки из первой части работы. Кроме того, по второй части предложения о подклееном пути из первой части работы, 
 оставшаяся часть пути не более $\frac{1}{3}(2^{s+2})$. Но часть пути после $Z_3$ не меньше, чем в четыре раза длиннее. Получаем противоречие, значит слово можно привести к нулю.


Пусть $V$ не содержит скобок вообще. В этом случае мы имеем дело с плоской ситуацией, без выходов в подклейки, также нет и входов и выходов в подплитки.




В этом случае, все ребра главные. Будем следить, по каким ребрам проходит предполагаемый путь. Допустим, $V$ содержит узловую букву кодирующую не боковой узел (то есть, узел типа $\mathbb{A}$, $\mathbb{B}$, $\mathbb{C}$, или угловой), тогда путь должен проходить по нескольким ребрам одной макроплитки. Заметим, что нельзя посетить более четырех внутренних ребер макроплитки не повернув ни в одном месте назад и не выйдя на периметр. Но $V^8$ должен содержать больше четырех таких участков, значит, в каком-то из них можно будет применить мономиальное соотношение.

Пусть теперь все узловые буквы $V$ кодируют боковые узлы. То есть предполагаемый путь проходит по одному и тому же ребру некоторой макроплитки. Для такого случая можно применить предложение~\ref{checkedges}. То есть мы применяем несколько преобразований, и если в процессе нельзя было применить мономиальное соотношение, в конце получится представление нашего слова в виде кодировки некоторого пути на комплексе. То есть, наш путь лежит полностью внутри некоторого ребра на комплексе.

\medskip

Внутри каждого ребра на комплексе нет более трех узлов с типами, отличными от $\mathbb{UL}$ и $\mathbb{LU}$.  Путь $V^8$ содержит восемь одинаковых участков подряд, в них не может встретиться вершин с типами, отличными от $\mathbb{UL}$ и $\mathbb{LU}$.
Но тогда в слове не может быть даже двух одинаковых участков подряд (лемма~\ref{nosquares}).
Таким образом, $W^9$ всегда приводится к нулю.

\end{proof}

Теперь мы можем завершить доказательство нашей основной теоремы:

\begin{theorem}
В полугруппе $S$, заданной конечным числом определяющих соотношений, существует бесконечное число различных слов, не равных нулю. При этом для любого слова $W$ его девятая степень приводится к нулю: $W^9=0$.

\end{theorem}

\begin{proof} Каждому регулярному пути на комплексе можно привести в соответствие его кодировку -- слово в полугруппе. При этом, при преобразовании пути его длина не меняется. Таким образом, если некоторый путь приводится к нулю, это значит, что либо в нем нашлось запрещенное подслово, либо локально некратчайший кусок пути, либо мертвый паттерн. Во всех этих случаях путь не является регулярным. То есть, все регулярные пути к нулю привести нельзя.

Вторая часть утверждения следует из леммы~\ref{powerpath}. Таким образом, построенная конечно определенная полугруппа $S$ содержит бесконечное множество различных слов, и при этом является нильполугруппой, где каждое слово в девятой степени приводится к нулю.

\end{proof}


{\bf Комментарий о вкладе авторов.}

Первому автору (Иванов-Погодаев) принадлежит конструкция подстановочной системы семейства комплексов, система кодировки вершин, перебор для проверки детерминированности и  финальная схема приведения периодического слова к нулю.

\section{Дальнейшие вопросы и соображения} \label{future}

Методы данной работы, возможно, помогут построить и иные объекты. В частности, представляет интерес следующий

\medskip {\bf Вопрос.} {\it Существует ли  конечно определенная полугруппа с размерностью Гель\-фан\-да--Ки\-рил\-ло\-ва $2.5$.}

\smallskip

Ранее были уже получены такие примеры для размерности $10,5$, \cite{BI}, \cite{ivadis}. Ранее для таких примеров применялась автоматная техника, которая также позволяла решить некоторые вопросы: \cite{ivamalev}, \cite{ivamalevsapir}.

\medskip

Л.~Н.~Шеврин поставил следующие вопросы:

\medskip
{\bf Вопрос 1.}
{\it Каково наименьшее возможное число элементов базиса бесконечной конечно определенной нильполугруппы конечного индекса?}

\medskip
{\bf Вопрос 2.} {\it Каков наименьший возможный индекс бесконечной конечно определенной нильполугруппы конечного индекса?}
\medskip

 Можно ли опустить в формулировке обоих вопросов слова «конечного индекса»? А именно, эти слова нужно будет опустить, если будет дан положительный ответ на следующий вопрос:

\medskip
{\bf Вопрос 3.} {\it Любая ли бесконечная конечно определенная полугруппа имеет конечный индекс?}

Нам представляются чрезвычайно важным для дальнейшего обобщение  результата Х.~Гудмана-Штраусса. В частности, представляют интерес два вопроса:

\medskip
{\bf А.} Топологическое обобщение результата Х.~Гудмана-Штраусса. Мозаика рассматривается как плоский граф, где плитками являются простые циклы. Нужно доказать, что заданием конечного числа локальных правил можно добиться заданного иерархического правила построения.

\medskip
{\bf В.} Обобщение иерархического правила. Нужно доказать, что заданием конечного числа локальных правил можно добиться заданного рекурсивного правила построения мозаики.

\medskip

Последнее нам представляется важным как в связи с третьим вопросом Л.~Н.~Шеврина (мы предполагаем отрицательный ответ), так и для построениия иных рекурсивных объектов. По всей видимости, эти вопросы должны предшествовать теоретико-кольцевым исследованиям.

И, наконец, возникает вопрос, на который также обратил наше внимание Л.~Н.~Шеврин:

\medskip
{\bf Вопрос 4.}  {\it Какие максимальные ненулевые степени элементов возможны в построенной нами полугруппе~$S$?}
\medskip

Все перечисленные вопросы представляются полезными для получения идей, позволяющих работать в кольцевой тематике.






\end{document}